# ESTIMATION OF TREND IN STATE-SPACE MODELS: ASYMPTOTIC MEAN SQUARE ERROR AND RATE OF CONVERGENCE

By Prabir Burman and Robert H. Shumway

*University of California, Davis*

The focus of this paper is on trend estimation for a general state-space model $Y_t = \mu_t + \varepsilon_t$, where the $d$th difference of the trend $\{\mu_t\}$ is assumed to be i.i.d., and the error sequence $\{\varepsilon_t\}$ is assumed to be a mean zero stationary process. A fairly precise asymptotic expression of the mean square error is derived for the estimator obtained by penalizing the $d$th order differences. Optimal rate of convergence is obtained, and it is shown to be "asymptotically equivalent" to a nonparametric estimator of a fixed trend model of smoothness of order $d - 0.5$. The results of this paper show that the optimal rate of convergence for the stochastic and nonstochastic cases are different. A criterion for selecting the penalty parameter and degree of difference $d$ is given, along with an application to the global temperature data, which shows that a longer term history has nonlinearities that are important to take into consideration.

**1. Introduction.** Trend estimation for time series data has a long history, and the literature, understandably, is quite vast. The basic statistical model and the estimation problem are quite easy to describe. The observed series $\{Y_t : t = 1, \ldots, n\}$ is modeled as

$$Y_t = \mu_t + \varepsilon_t, \tag{1.1}$$

where the error series $\{\varepsilon_t\}$ is assumed to be a mean zero stationary process. In some cases, the error series may have variances that change with time, but we will not worry about that issue here. The goal is to estimate the trend $\{\mu_t\}$ on the basis of the observed data. In the literature, two types of structures of the trend are assumed: fixed and stochastic. Asymptotic analysis of the estimator of the random trend model (a version of the state-space model) is the focus of this paper. We derive the expression for the asymptotic mean square error of the trend estimate obtained by penalizing finite differences and then obtain its rate of convergence. The main asymptotic results presented in this paper are not in terms of upper bounds, rather they are asymptotic expressions, which are correct up to first order with bounds on the second order terms.

In the fixed trend model, it is usually assumed that the trend is of the form $\mu_t = \mu(t/n)$, where $\mu$ is an unknown function (presumably smooth) on the interval $[0,1]$. Various methods, such as kernel [Altman (1990) and Truong (1991)], local polynomial [Beran and Feng (2001)], spline [Burman (1991)] and wavelets [Johnstone and Silverman (1997)], have been employed in order to estimate the trend. Asymptotic properties of such methods, rates of convergence and issues on smoothing parameter selection have been well investigated by many authors [see, e.g., Eubank (1988), Fan and Yao (2003), Tran et al. (1996) and Robinson (1997)]. The literature is vast, and the recent book by Fan and Yao is a good source on this topic and the relevant references.

The stochastic trend models are quite popular in the time series literature [see Chapter 9 in Box, Jenkins and Reinsel (1994), Chapter 3 in Durbin and Koopman (2001), Harvey (1991) and Chapter 4 in Shumway and Stoffer (2000)]. While deterministic trends may sometimes have simpler expressions, it is the belief of many time series analysts that a stochastic trend model is more realistic in real applications. A discussion on deterministic versus stochastic trend can be found in Chapter 4.1 in the book by Box, Jenkins and Reinsel (1994). For the random trend model, it is assumed that the $d$th order differences

$$\nabla^d \mu_t = \gamma_t \tag{1.2a}$$

are mean zero i.i.d. variables with variance $\sigma_\gamma^2$, where $d > 0$ may or may not be an integer. The usual application in the literature assumes that $d$ is a known positive integer. Independence of $\{\nabla^d \mu_t\}$ is not really necessary, stationarity is enough (see Remark 1 in Section 2). In this paper, we restrict our attention to the case when $d$ is an integer. In a forthcoming paper, we will deal with the general case $d > 0$ in detail. An alternative representation for the random trend is

$$\mu_t = \sum_{0 \le j \le d-1} \beta_j t^{j-1} + \sum_{1 \le j \le t} (-1)^{t-j} \binom{-d}{t-j} \gamma_j, \tag{1.2b}$$

where $\beta_j$'s, the coefficients of the polynomial, can be taken to be fixed or random. Chan and Palma (1998) have investigated finite approximations to the log-likelihood for state-space models in the long memory (fractional difference) case. State-space models have been quite popular among many statisticians, and efficient algorithms such as EM and MCMC have been developed for estimation of the trend [see, e.g., Shumway and Stoffer (2000) and Durbin and Koopman (2001)]. However, not much is known about the asymptotic properties of the estimators, even though there is a strong parallel with the method of smoothing splines [see, e.g., Wahba (1978), Wecker and Ansley (1983) and Kohn, Ansley and Wong (1992)].

Assuming that $\{\varepsilon_t\}$ are i.i.d. $N(0, \sigma_\varepsilon^2)$ and $\{\nabla^d \mu_t\}$ are i.i.d. $N(0, \sigma_\gamma^2)$, we note that (1.1) and (1.2a) define a special case of the state space model, where (1.1) is the observation equation and (1.2) is the state equation. If the errors are not i.i.d. but stationary, a second state equation defining $\{\varepsilon_t\}$ as a stationary autoregressive process can be added. For the i.i.d. error case, the estimator of the trend is obtained by maximizing the Gaussian likelihood; that is, by minimizing

$$(1.3) \qquad \sum_{1 \le t \le n} (Y_t - \mu_t)^2 + \nu \sum_{d+1 \le t \le n} (\nabla^d \mu_t)^2$$

with respect to $\mu = (\mu_1, \ldots, \mu_n)'$, where $\nu = \sigma_\varepsilon^2 / \sigma_\gamma^2$. The estimate, for fixed $\nu$ and $d$, is given by $\hat{\mu}_t = E(\mu_t | Y_1, \ldots, Y_n)$, which can also be computed as the Kalman smoothers. Note that, restricting the class to linear estimators, the above minimization problem is still valid without the Gaussian assumptions on $\{\varepsilon_t\}$ and $\{\gamma_t\}$, and the resulting estimator in such a case can be described as the best linear unbiased predictor (BLUP) of $\{\mu_t\}$ [Kimeldorf and Wahba (1970)].

Wahba (1978) established a connection to spline smoothing. Suppose that the trend is given by

$$(1.4) \qquad \mu_t = \sum_{0 \le j \le d-1} \theta_j t^j + \tau \int_0^t (t-u)^{d-1} dW(u),$$

where $W$ is a standard Brownian motion and $\theta = (\theta_0, \ldots, \theta_{d-1})'$ is multivariate normal $N_d(a, kI)$. Assuming a diffuse prior for $\theta$, that is, if $k \to \infty$, then minimizing

$$(1.5) \qquad \sum_{1 \le t \le n} (Y_t - \mu_t)^2 + \nu \int_0^t \left\{ \frac{\partial^d \mu_u}{\partial u^d} \right\}^2 du$$

with respect to $\mu_1, \ldots, \mu_n$ leads to the estimator $\hat{\mu}_t = E(\mu_t | Y_1, \ldots, Y_n)$, which is precisely the same estimator as one obtains for the state-space model [see Wahba (1978) and Green and Silverman (1994) for details]. In other words, under the diffuse prior model of (1.4), the smoothing spline estimate of $\mu_t$ obtained by minimizing (1.5) coincides with the estimator obtained by a criterion that involves penalizing finite differences. However, the properties of the latter estimator are still not well understood.

As mentioned in the first paragraph, this paper is devoted to the asymptotic study of the estimator $\hat{\mu}$ obtained by minimizing (1.3) (and its weighted least squares variant). Specifically, we obtain an asymptotic expression of the mean square error under fairly general assumptions. We only assume that $\{\varepsilon_t\}$ is mean zero stationary whose spectral density at zero is positive and that $\{\nabla^d \mu_t\}$ are i.i.d. Note that no assumption of Gaussianity is needed for our results to hold. We show that the smallest mean square error is constant times $n^{2(d-0.5)/(2d)}$. This rate of convergence is the same as one obtains for the fixed trend case when the order of smoothness is $d-0.5$. It may be worthwhile to point out that the rates for the stochastic and nonstochastic cass are different. For the nonstochastic case, the rate of convergence would be of order $n^{2d/(2d+1)}$. The penalized least squares method proposed here does not require the knowledge of whether the underlying trend is deterministic or stochastic. So, if the underlying trend is stochastic, but one mistakes it be deterministic and applies the method proposed here, all the results given in this paper remain valid. However, if the underlying trend is deterministic, then it is reasonable to believe (following the theory of nonparametric function estimation) that the optimal rate associated with our procedure is of order $n^{2d/(2d+1)}$.

In Section 2, we state the main results where the asymptotic expressions for the mean square error of the trend estimates for the state-space model are obtained. We discuss two types of estimates: ordinary least squares and the weighted least squares. In Section 3, we present a criterion for estimating the penalty parameter $\nu$ and order of difference $d$ along with a numerical example. Finally in the Appendix we present the proofs of the results of Section 2. The proofs of the results depend on the properties of Toeplitz, Hankel and circulant matrices.

**2. The main results.** We now take up the issue of trend estimation in a state-space model. Two methods will be discussed: first, the ordinary least squares and, second, in Section 2.1, the weighted least squares. Clearly, these two methods yield the same estimate when the error sequence is assumed to be i.i.d. Consider the time series $Y_t = \mu_t + \varepsilon_t$, $t = 1, \ldots, n$, where $\mu_t$ is the trend and $\varepsilon_t$ is a mean zero stationary process. We will assume that the $d$th difference of the trend is i.i.d. In the literature, it is usually assumed that the errors $\varepsilon_t$ are i.i.d. Gaussian and the $d$th differences of the trend $\nabla^d \mu_t = \gamma_t$ are also i.i.d. Gaussian variables. In such a case, the trend can be estimated by minimizing the negative of the log-likelihood; that is, by minimizing (1.3) with respect to $\mu_t$'s, $\sigma_\varepsilon$ and $\sigma_\gamma$. There are a few points to be noted here. The assumption of normality can be dispensed with, and, in such a case, we can treat this as a problem of obtaining the best prediction of linear predictor (BLUP) of the random effects $\gamma_t$'s for the mixed linear model. Throughout this paper, the error sequence $\{\varepsilon_t\}$ is assumed to be stationary. So, instead of minimizing the penalized ordinary least squares criterion given in (1.3), we should perhaps minimize a weighted least squares criterion, and this is taken up in Section 2.1. In this subsection, we still deal with the ordinary least squares criterion given in (1.3) in the presence of stationary error series $\{\varepsilon_t\}$. In the next section, we will describe a method for estimating $\nu$ (and the degree of difference $d$).

Before we begin, it is important to note that we would like the average signal-to-noise ratio $E\|\mu\|^2/(n\sigma_\varepsilon^2)$ to be bounded between two positive constants. This can be guaranteed if the quantity $E\|\mu\|^2/n$ stays bounded between two positive constants. This requirement is met if the $d$th order difference $\gamma_t = \nabla^d \mu_t$ has a variance of the form

$$\sigma_\gamma^2 = \tau^2/n^{2d-1} \tag{2.1}$$

for some positive constant $\tau$. In order to see why this is true, we need to rewrite the random trend in (1.2b) as

$$\mu_t = \sum_{0 \le j \le d-1} \beta_j (t/n)^{j-1} + \sum_{1 \le j \le t} (-1)^{t-j} \binom{-d}{t-j} \gamma_j = \mu_{1t} + \mu_{2t} \qquad \text{say.}$$

First, let us consider the polynomial part. If the coefficients $\beta_j$'s are norandom and are bounded in absolute values by a constant independent of $n$, then $\|\mu_1\|^2/n$ is finite. If $\beta_j$'s are random with finite means and variances (which do not depend on $n$), then $E\|\mu_1\|^2/n$ is also finite. Now, let us consider the purely random part $\mu_{2t}$ of the trend. By Stirling's approximation, $(-1)^l \binom{-d}{l}$ is approximately equal to $(\Gamma(d))^{-1} l^{d-1}$ for a positive integer $l$ not small, where $\Gamma$ is the usual Gamma function. Then, a fairly straightforward calculation will show that $E\|\mu_2\|^2/n$ is approximately equal to constant times $n^{2d-1}\sigma_\gamma^2$, when $n$ is not small.

The statistical model, which we assume to be correct throughout, is

$$\begin{aligned}
& Y_t = \mu_t + \varepsilon_t \qquad \text{with } \nabla^d \mu_t = \gamma_t, \text{ where} \\
& \{\varepsilon_t\} \text{ is mean zero stationary,} \\
& \{\gamma_t\} \text{ are i.i.d. with mean zero and variance } \tau^2/n^{2d-1}, \\
& \{\varepsilon_t\} \text{ and } \{\gamma_t\} \text{ are independent.}
\end{aligned} \tag{2.2}$$

We have discussed, above, that the trend consists of a polynomial part $\mu_{1t}$ and a purely random part $\mu_{2t}$. It is important to point out that the purely random part $\mu_{2t}$ is a nonnegligible one. This can be seen once we note that $\mu_{2t}$ has zero mean and $E(\mu_{2t}^2)$ approximately equal to a constant times $(t/n)^{2d-1}$. Thus, the purely random part of the trend is nonnegligible except when $t$ is small.

We will now find a matrix representation of the estimate of $\mu$ obtained by minimizing (1.3). For $d > 0$, the summation and difference operators on

$R^n$ denoted by $S_d$ and $S_{-d}$, respectively, are defined as follows: for any $x$ in $R^n$,

$$(S_d x)(t) = \sum_{1 \leq j \leq t} (-1)^{t-j} \binom{-d}{t-j} x_j, (S_{-d}x)(t) \tag{2.3}$$

$$= \sum_{1 \leq j \leq t} (-1)^{t-j} \binom{d}{t-j} x_j.$$

It can be shown that $S_d S_{-d} = I$ and $S_0 = I$. Properties of summation and difference operators can be found in Burman (2006). Note that the difference operator $S_{-d}$ is a lower triangular band matrix with element $(t,j)$ is $(-1)^{t-j}\binom{d}{t-j}$, $t \geq j$. Let the $(n-d) \times n$ matrix obtained by deleting the first $d$ rows of $S_{-d}$ be denoted by $\overline{S}_{-d}$ as follows. Then, with $Y = (Y_1, \ldots, Y_n)'$, we can rewrite (1.3) as

$$SSE(\mu, \nu, d) = \|Y - \mu\|^2 + \nu \mu' \overline{S}'_{-d} \overline{S}_{-d} \mu. \tag{2.4}$$

Hence, minimizing the expression in (2.4) with respect to $\mu_t$'s leads to the estimate

$$\hat{\mu} = \hat{\mu}(\nu) = (I + \nu \overline{S}'_{-d} \overline{S}_{-d})^{-1} Y. \tag{2.5}$$

The estimate $\hat{\mu}$ given above in (2.5) is rather easy to calculate, even for large $n$, since $I + \nu \overline{S}'_{-d} \overline{S}_{-d}$ is a band matrix. It should also be emphasized that this matrix representation eliminates the need for deciding the intitial values for the Kalman iterations.

It is possible to construct a pointwise prediction interval for $\mu_t$ under Gaussian assumption when the errors $\{\varepsilon_t\}$ are assumed to be i.i.d., and the variances $\sigma_\varepsilon^2$ and $\sigma_\gamma^2$ are known. In such a case, we can express the model as $Y = X\beta + Z\gamma + \varepsilon$, where element $(t,j)$ of the matrix $X$ is $(t/n)^{j-1}$, $j = 1, \ldots, d$, $Z$ is a lower triangular matrix of order $n$ whose element $(t,j)$ is $(-1)^{t-j}\binom{-d}{t-j}$. The conditional distribution of $\mu$ given $Y$ is normal with mean

$$E(\mu|Y) = X\beta + (I + \nu^*(Z'Z)^{-1})^{-1}(Y - X\beta)$$

and variance–covariance matrix

$$D = \sigma_\varepsilon^2 (I + \nu^*(ZZ')^{-1})^{-1},$$

where $\nu^* = \sigma_\varepsilon^2/\sigma_\gamma^2$. It can be shown that the estimate given in (2.5) is $\hat{\mu}(\nu) = X\hat{\beta} + (I + \nu^*(Z'Z)^{-1})^{-1}(Y - X\hat{\beta})$, where $\hat{\beta}$ is the estimate of $\beta$ obtained by the weighted least squares criterion $(Y - X\beta)'(\sigma_\varepsilon^2 I + \sigma_\gamma^2 ZZ')^{-1}(Y - X\beta)$. Since $\hat{\beta}$ is a $\sqrt{n}$-consistent estimate of $\beta$, the conditional mean of $\mu$ given $Y$ is well estimated by $\hat{\mu}(\nu)$. An approximate $(1-\alpha)100\%$ prediction interval of $\mu_t$ is given by $\hat{\mu}_t \pm z_\alpha D_{tt}$, where $z_\alpha$ is the critical value from the standard normal distribution and $\{D_{tt}\}$ are the diagonal elements of the variance–covariance $D$ matrix of the conditional distribution of $\mu$ given $Y$. In practice, however, $\sigma_\varepsilon$ and $\nu^* = \sigma_\varepsilon^2/\sigma_\gamma^2$ are unknown and have to be estimated from the data. In Section 3, we will discuss these issues. The matrix $D$ is not as formidable as it looks. It is a banded matrix, since $Z^{-1}$ is a banded lower triangular matrix whose element $(t,j)$ is given by $(-1)^{t-j}\binom{d}{t-j}$, $0 \leq t-j \leq d$.

It is not unusual to have time series data with missing observations. If a large block of consecutive observations are missing, then nothing can be done to estimate the trend during those periods of time. Suppose, for the case of simplicity, that all the observations between time periods $n_1$ and $n_2$ are missing, where $n_1 = np_1$ and $n_2 = np_2$, and $0 < p_1 < p_2 < 1$. Then, the methods of this paper can be used to estimate the trend $\mu_t$ for $1 \leq t \leq n_1$ and $n_2 \leq t \leq n$, and all the results will be valid with appropriate modifications. However, the more difficult case is when the obervations are sporadically missing. In such a case, the entire trend $\{\mu_t : t = 1, \ldots, n\}$ should be estimable. A reasonable approach for estimating the trend would be to minimize

$$\sum_{t \in J}(Y_t - \mu_t)^2 + \nu \mu' \overline{S}'_{-d} \overline{S}_{-d} \mu$$

with respect to $\mu$ where $J$ is set of time indices at which the observations are available. In such a case, the estimate of $\mu$ would be $\hat{\mu} = (\widetilde{I} + \nu \overline{S}'_{-d} \overline{S}_{-d})^{-1} Y$, where $\widetilde{I}$ is a diagonal matrix whose $t$th diagonal element is 1 or 0 depending on whether the observation is available or missing. We have not yet investigated this estimate and its properties, and this issue needs further research.

2.1. *Asymptotic mean square error of the estimate of the random trend.* Let us denote the conditional mean of the estimate, given the trend as $\overline{\mu}_t = E(\hat{\mu}_t | \mu_1, \ldots, \mu_n)$. If we view $\overline{\mu}_t - \mu_t$ as the bias of the estimate $\hat{\mu}_t$, then we can decompose the mean the square error as in the usual variance-bias decomposition

$$E\|\hat{\mu} - \mu\|^2 = E\|\hat{\mu} - \overline{\mu}\|^2 + E\|\overline{\mu} - \mu\|^2.$$

Theorems given below tell us the asymptotic values of the bias and the variance. We will provide an outline of the results here. If we write $b = \nu^{1/(2d)}/n$ and assume that $b \to 0$ and $nb \to \infty$ as $n \to \infty$, then it turns out that

$$E\|\hat{\mu} - \overline{\mu}\|^2/n = c_1(nb)^{-1}[1 + O((nb)^{-1}) + O(b)],$$
$$E\|\overline{\mu} - \mu\|^2/n = c_2 b^{2d-1}[1 + O((nb)^{-1}) + O(b)],$$

where $c_1$ and $c_2$ are constants, the expressions for which are to be found in Theorems 1 and 2 given below. Note that the results look very much like the asymptotic expressions of the variance and bias-square components for estimating a regression function using a kernel method with bandwidth $b$. Also, note that the bias looks like that of a function (nonrandom), which is $d-1$ times differentiable with the $(d-1)$st derivative satisfying a Lipschitz condition of order 0.5. Hence, we can obtain the value of $b = b^*$ at which the mean square error is minimized and calculate the minimum mean square error explicitly. More discussion of the results are given below after the theorems have been stated.

The first result given below is on the asymptotic expression of the variance. The proofs of both the theorems need to employ the theories of Toeplitz, Hankel and circulant matrices.

THEOREM 1. *Assume that the conditions given in (2.2) hold and that $\sum_{1 \leq j < \infty} j|\rho_\varepsilon(j)| < \infty$, where $\rho_\varepsilon(j)$ is the covariance of the stationary process of lag $j$ for the error process $\{\varepsilon_t\}$. Assume that $g_\varepsilon(0) > 0$, where $g_\varepsilon(u) = \sum_{-\infty < j < \infty} \rho_\varepsilon(j) e^{iju}$ is the spectral density function of the error process. Then, assuming $\nu \to \infty$ and $\nu/n^{2d} \to 0$ as $n \to \infty$, we have*

$$E\|\hat{\mu} - \overline{\mu}\|^2/n = c_1 \nu^{-1/(2d)}[1 + O(\nu^{-1/(2d)}) + O(\nu^{1/(2d)}/n)],$$

*where $c_1 = g_\varepsilon(0) \operatorname{Beta}(1/(2d), 2 - 1/(2d))/(2d\pi)$.*

The following result gives an asymptotic expression for the bias.

THEOREM 2. *Assume that the conditions given in (2.2) hold. Assuming that $\nu \to \infty$ and $\nu/n^{2d} \to 0$ as $n \to \infty$, we have*

$$E\|\overline{\mu} - \mu\|^2/n = c_2 (\nu^{1/(2d)}/n)^{2d-1}[1 + O(\nu^{-1/(2d)}) + O(\nu^{1/(2d)}/n)],$$

*where $c_2 = \tau^2 \operatorname{Beta}(1 + 1/(2d), 1 - 1/(2d))/(2d\pi)$.*

REMARK 1. We have so far assumed that the $d$th difference $\{\gamma_t = \nabla^d \mu_t\}$ of the trend consists of i.i.d. random variables with mean zero and variance $\sigma_\gamma^2 = \tau^2/n^{2d-1}$. Are Theorems 1 and 2 valid when $\{\nabla^d \mu_t\}$ are not i.i.d.? The answer is yes, if we assume that $\{\gamma_t = \nabla^d \mu_t\}$ is a mean zero stationary process with autocovariances $\operatorname{Cov}(\gamma_{t+j}, \gamma_t) = \sigma_\gamma^2 \rho_\gamma(j)$, $-\infty < j < \infty$. In such a case, Theorem 1 is exactly the same as before. Theorem 2 is also the same as before except for the constant that appears in the expression of $E\|\overline{\mu} - \mu\|^2/n$. Let $g_\gamma$ be the spectral density of the process $\{\gamma_t/\sigma_\gamma\}$ and assume that $\sum_{1 \leq j < \infty} j|\rho_\gamma(j)| < \infty$. A modification of the proof of Theorem 2 will show

$$E\|\overline{\mu} - \mu\|^2/n = c_2 g_\gamma(0)(\nu^{1/(2d)}/n)^{2d-1}[1 + O(\nu^{-1/(2d)}) + O(\nu^{1/(2d)}/n)],$$

where $c_2$ is the same constant as in Theorem 2. Note that the constant involved in the expression of $E\|\overline{\mu} - \mu\|^2/n$ now includes the value of the spectral density $g_\gamma$ at zero. Clearly, for the case when $\{\nabla^d \mu_t\}$ is stationary, all the discussion below about the optimal mean square error in estimating the trend remains valid with appropriate constants.

REMARK 2. Note that the mean square error
$$D(\nu) = E\|\hat{\mu} - \mu\|^2/n = [c_1 \nu^{-1/(2d)} + c_2 \nu^{1-1/(2d)} n^{-2d+1}](1 + o(1))$$
is minimized at
$$\nu^* = n^{2d-1}(2d-1)^{-1}(c_1/c_2)(1 + o(1))$$
and the smallest mean square error is
$$D(\nu^*) = n^{(2d-1)/(2d)} c_3 (1 + o(1)),$$
where $c_3 = c_1(c_2/c_1)^{1/(2d)} 2d(2d-1)^{-1+1/(2d)}$.

REMARK 3. We consider, here, the Euclidean distance between the true random trend $\mu$ and its estimate $\hat{\mu}$. Such a distance has been considered by many for time-dependent observations [see, e.g., Altman (1990), Burman (1991), Johnstone and Silverman (1997) and Truong (1991)]. However, it is of interest to consider the distance $(\hat{\mu} - \mu)' R_\varepsilon^{-1} (\hat{\mu} - \mu)$, where $R_\varepsilon$ is the $n \times n$ variance–covariance matrix of the error process $\{\varepsilon_t\}$. If the spectral density function of this process stays bounded between two positive constants, which is the case for the usual ARMA model, then the theory of Toeplitz matrices tells us that all the eigenvalues of the matrix $R_\varepsilon$ stay between two positive constants [Grenander and Szego (1958)]. In such a case, we can find two positive constants $k_1$ and $k_2$ such that
$$k_1 \|\hat{\mu} - \mu\|^2 \leq (\hat{\mu} - \mu)' R_\varepsilon^{-1} (\hat{\mu} - \mu) \leq k_2 \|\hat{\mu} - \mu\|^2.$$
Consequently, all of the rates of convergence results for the Euclidean distance $\|\hat{\mu} - \mu\|^2$ are also valid when the distance is taken to be $(\hat{\mu} - \mu)' R_\varepsilon^{-1} (\hat{\mu} - \mu)$.

We will conclude this subsection by comparing the optimal mean square error as discussed above with the optimal rate of convergence associated with nonparametric function estimation problems. A function $f$ on $[0,1]$ is defined to be in the smoothness class $p = r + \beta$, where $r$ is a nonnegative integer and $0 < \beta \leq 1$, if $f$ is $r$ times differentiable and the $r$th derivative of $f$ is Lipschitz of order $\beta$. Now, if the trend is a nonrandom function and is modeled as $\mu_t = \mu(t/n)$, and the function $\mu$ is of smoothness class $p$, then the optimal rate of convergence for estimating the trend is given by $n^{-2p/(2p+1)}$ [see Stone (1982), Eubank (1988) and Fan and Yao (2003)]. Theorems 1 and 2 and the subsequent discussion tell us that, for the state-space model as given in (2.2), the optimal rate of convergence is $n^{-(2d-1)/(2d)}$. This corresponds to the rate of $p = d - 0.5$ for the nonrandom case.

*What is the rate when the order of difference in unknown?* We now examine the performance of the estimated stochastic trend when the true order of difference, which we assume to be $d_0$, is unknown, but a $d$th order differencing scheme is being employed to estimate the trend. It may be worthwhile to point out that the order of difference $d$ controls the smoothness of the estimated trend. Since, in the time series literature, it is assumed that the true order of difference is known, we may turn to the literature on nonparametric function estimation for some guidance. For nonparametric function estimation, it has been pointed out by many authors that, from a practical point of view, the penalty parameter (or the smoothing parameter in general) is far more important than the order of differencing $d$ employed in the estimation procedure [see Beran (2005), Eubank (1988) and Wahba (1990)]. For instance, in his analysis of multi-way tables, Beran finds that there is not much of a difference in the estimated risk when the order of difference is taken to be any integer between 1 and 4.

However, for the sake of theoretical completeness, we will present the results when the true difference order $d_0$ is unknown, and we are employing a $d$th order differencing in order to obtain the estimate of the stochastic trend. Before we state the result, let us point out that, when $d \geq d_0$, the rate remains the same though the constant associated with the rate depends on $d$. However, the constant is the smallest when $d = d_0$. When $d < d_0$, a different rate comes into play, and this rate is the same as in the usual nonparametric function estimation [see, e.g., Eubank (1988), Stone (1982) and Wahba (1990)].

It should be pointed out that Theorem 1 remains valid even when $d \neq d_0$. However, Theorem 2 is no longer valid when $d \neq d_0$. We will write the results for the bias part (i.e., analogue of Theorem 2). When $d < d_0$, we can only obtain a bound for the bias part. But, for the case $d \geq d_0$, we can obtain an asymptotic expresssion.

THEOREM 3. *Assume that the model as given in (2.2) holds with $d$ replaced by $d_0$. The estimate of $\mu$ is obtained by minimizing the expression given in (1.3):*

(a) *When $d < d_0$,*

$$E\{\|\overline{\mu} - \mu\|^2\}/n = O((\nu/n)^{2d}).$$

(b) *When $d \geq d_0$,*

$$E\{\|\overline{\mu} - \mu\|^2\}/n = c_4(\nu^{1/(2d)}/n)^{2d_0-1}[1+o(1)],$$

*where $c_4 = \tau^2 \operatorname{Beta}((2d_0 - 1)/(2d), 2 - (2d_0 - 1)/(2d))/(2d\pi)$.*

Even though we obtain an upper bound for the case $d < d_0$, we believe it is not possible to improve the rate, and that is certainly true when $d = 1$, the case for which the exact expressions for the eigenvalues and eigenvectors for the matrix $\overline{S}'_{-d}\overline{S}_{-d}$ are known. Note that the optimal mean square error $E\{\|\hat{\mu} - \mu\|^2\}/n$ is $O(n^{-2d/(2d+1)})(1 + o(1))$, which is known to be the rate for nonparametric function estimation for the nonstochastic case.

When $d \geq d_0$, Theorem 1 and part (b) of Theorem 3 tell us that the mean square error is of the form $c(d)n^{-(2d_0-1)/(2d_0)}$ when

$$c(d) = [g_\varepsilon(0)^{2d_0-1}\tau^2(2d_0-1)(\operatorname{Beta}(1/(2d) - 2 - 1/(2d))^{2d_0-1}$$
$$\times \operatorname{Beta}((2d_0-1)/(2d), 2 - (2d_0-1)/(2d))]^{1/(2d_0)}.$$

It can be shown that $c(d)$ is minimized when $d = d_0$. In other words, the optimal rate of covergence is still $n^{-(2d_0-1)/(2d_0)}$ as long as $d \geq d_0$, but the constant associated with the rate depends on $d$ and the minimum value of the constant is acheived at $d = d_0$.

2.2. *Weighted least squares estimate of the trend.* In this subsection, we will discuss a weighted least squares estimator of the trend and its asymptotic properties. Since the arguments needed to prove the results given in this subsection are similar to those for Theorems 1 and 2, we will merely state the results. We will first discuss the case when the variance–covariance matrix $R_\varepsilon$ of $\{\varepsilon_t : t = 1, \ldots, n\}$ is assumed to be known. A weighted least squares estimate $\widetilde{\mu}^{(\mathrm{wls})}$ of $\mu$ may be obtained by minimizing

$$(2.6) \qquad (Y - \mu)' R_\varepsilon^{-1}(Y - \mu) + \nu \sum_{d+1 \leq t \leq n}(\nabla^d \mu_t)^2$$

instead of minimizing the quantity given in (1.3). Clearly, then

$$(2.7) \qquad \widetilde{\mu}^{(\mathrm{wls})} = (R_\varepsilon^{-1} + \nu \overline{S}'_{-d}\overline{S}_{-d})^{-1} R_\varepsilon^{-1} Y.$$

In practice, of course, the matrix $R_\varepsilon^{-1}$ is unknown and has to be estimated from data. If $\hat{R}_\varepsilon$ is an estimate of $R_\varepsilon$, then a practical weighted least squares estimate of the trend is given by replacing $R_\varepsilon$ in (2.7) by $\hat{R}_\varepsilon$, and we denote the resulting estimator by $\hat{\mu}^{(\mathrm{wls})}$. We will obtain, below, analogues of Theorems 1 and 2 for $\widetilde{\mu}^{(\mathrm{wls})}$ and also show that the difference $\hat{\mu}^{(\mathrm{wls})} - \widetilde{\mu}^{(\mathrm{wls})}$ is small in the probabilistic sense under appropriate conditions.

We will assume that the error process $\{\varepsilon_t\}$ is $\mathrm{AR}(p)$ [or $\mathrm{ARMA}(p, q)$]. Using a preliminary estimate $\hat{\mu}^+$ of $\mu$, we can estimate $p$ (via a model selection criterion such as AIC or BIC) and the parameters of the error process by using the residuals $Y_t - \hat{\mu}_t^+$. Clearly, an estimate of the variance–covariance matrix $R_\varepsilon$ can be obtained from the estimated model of the error process. Moreover, the estimated variance–covariance matrix $\hat{R}_\varepsilon$ is a $\sqrt{n}$ consistent estimate of $R_\varepsilon$; that is, $\|\hat{R}_\varepsilon - R_\varepsilon\| = O_P(n^{-1/2})$, where $\|\cdot\|$ is the usual matrix norm (i.e., the maximum of singular values).

For the following results on the asymptotics of the weighted least squares, for any $x$ in $R^n$, we define $\|x\|^2_{R_\varepsilon^{-1}}$ to be equal to $x'R_\varepsilon^{-1}x$.

THEOREM 4. *Assume that the conditions given in (2.2) holds and that the error process $\{\varepsilon_t\}$ is $\mathrm{AR}(p)$. Let $\overline{\mu}^{(\mathrm{wls})} = E[\widetilde{\mu}^{(\mathrm{wls})}|\mu]$. Then, as $\nu \to \infty$ and $\nu/n^{2d} \to 0$, the following results are true:*

(a) $\quad E[\|\widetilde{\mu}^{(\mathrm{wls})} - \overline{\mu}^{(\mathrm{wls})}\|^2_{R_\varepsilon^{-1}}]/n = c_4\nu^{-1/(2d)}[1 + O(\nu^{-1/(2d)}) + O(\nu^{1/(2d)}/n)],$

*where $c_4 = g_\varepsilon(0)^{-1/(2d)}\mathrm{Beta}(1/2d, 2 - 1/(2d))/(2d\pi)$;*

(b) $\quad E[\|\overline{\mu}^{(\mathrm{wls})} - \mu\|^2_{R_\varepsilon^{-1}}]/n$

$$= c_5\nu^{1-1/(2d)}n^{-2d+1}[1 + O(\nu^{-1/(2d)}) + O(\nu^{1/(2d)}/n)],$$

*where $c_5 = \tau^2 g_\varepsilon(0)^{-1/(2d)}\mathrm{Beta}(1 + 1/(2d), 1 - 1/(2d))/(2d\pi)$;*

(c) *The mean square error is given by*

$$E[\|\widetilde{\mu}^{(\mathrm{wls})} - \mu\|^2_{R_\varepsilon^{-1}}]/n$$
$$= [c_4\nu^{-1/(2d)} + c_5\nu^{1-1/(2d)}n^{-2d+1}][1 + O(\nu^{-1/d}) + O(\nu^{1/(2d)}/n)].$$

THEOREM 5. *Assume that the conditions for Theorem 3 hold. Moreover, assume that the fourth moment of the error process $\varepsilon_t$ is finite and $\|\hat{R}_\varepsilon - R_\varepsilon\| = O_P(n^{-1/2})$. Then,*

$$\|\widetilde{\mu}^{(\mathrm{wls})} - \mu\|^2_{R_\varepsilon^{-1}}/n - \|\hat{\mu}^{(\mathrm{wls})} - \mu\|^2_{R_\varepsilon^{-1}}/n = O_P(1/n).$$

REMARK 4. There are a number of consequences that follow from Theorems 3 and 4. First, from the asymptotic expression of the mean square error for the weighted least squares estimator $\widetilde{\mu}^{(\mathrm{wls})}$ of $\mu$ for a known $R_\varepsilon$, we can obtain the optimal rate of convergence. Note that the minimum mean square errors, for the weighted least squares estimate $\widetilde{\mu}^{(\mathrm{wls})}$ and the ordinary least squares estimate $\hat{\mu}$ given in Section 2.1, differ only in constants. Moreover, Theorem 4 guarantees that the weighted least squares estimate $\hat{\mu}^{(\mathrm{wls})}$ for unknown $R_\varepsilon$ has the same mean square error as $\widetilde{\mu}^{(\mathrm{wls})}$ in the asymptotic sense.

REMARK 5. Even though Theorem 3 is stated for the case when the error process $\{\varepsilon_t\}$ follows $\mathrm{AR}(p)$, this result is true for any stationary error process as long as its spectral density function is bounded away from 0 and $\infty$. Theorem 4 holds as long as the error process has a finite-dimensional model such as $\mathrm{AR}(p)$ or an invertible $\mathrm{ARMA}(p,q)$.

**3. Data dependent selection of $\nu$ and $d$.** In this section, we will discuss the issue of selecting the smoothing parameter $\nu$ and the degree of differences $d$ for the ordinary least squares estimate $\hat{\mu}$ given in Section 2.1. Criterion for selection of $\nu$ and $d$ can be developed for the weighted least squares estimate $\hat{\mu}^{(\text{wls})}$ by following similar arguments given in this section. However, we will not address that issue here. The main idea rests on minimizing the expected distance between $\mu$ and its estimate $\hat{\mu} = \hat{\mu}(\nu)$, given by

$$D(\nu, d) = E\|\hat{\mu} - \mu\|^2/n = E\|\hat{\mu} - \overline{\mu}\|^2/n + E\|\overline{\mu} - \mu\|^2/n,$$

where $\overline{\mu} = \overline{\mu}(\nu) = E[\hat{\mu}(\nu)|\mu]$. Ideally we would like to minimize $D$ with respect to $\nu$ (and $d$). However, it is unknown, and so we try the route of estimating $D$ by following the arguments given by Akaike and Mallows.

Now, the expected value of the residual sum of squares is

$$E(SSE(\nu, d)) = E\|Y - \hat{\mu}\|^2$$
$$= n\sigma^2 + E\|\hat{\mu} - \overline{\mu}\|^2 + E\|\overline{\mu} - \mu\|^2 - 2\operatorname{tr}((I + \nu\overline{S}'_{-d}\overline{S}_{-d})^{-1}R_\varepsilon),$$

where $\sigma^2 = E(\varepsilon_t^2)$ and variance–covariance matrix of the error series $\{\varepsilon_t\}$ is $R_\varepsilon$.

So, an unbiased estimate of $D(\nu, d)$ is given by

(3.1) $$\widetilde{D}(\nu, d) = SSE(\nu, d)/n + 2\operatorname{tr}((I + \nu\overline{S}'_{-d}\overline{S}_{-d})^{-1}R_\varepsilon)/n - \sigma^2.$$

Since the last term in the expression of $\widetilde{D}$ does not depend on $\nu$ (and $d$), we can safely ignore it. Using the same arguments as in Theorem 1, we can show that

$$\operatorname{tr}((I + \nu\overline{S}'_{-d}\overline{S}_{-d})^{-1}R_\varepsilon) = \sum g_\varepsilon(\pi j/n)/(1 + \nu s(\pi j/n)) + O(1),$$

where $g_\varepsilon$ is the spectral density of the process $\{\varepsilon_t\}$, $s(u) = (2 - 2\cos u)^d$, and

$$c_6(d) = \pi^{-1} \int_0^\infty 1/(1 + u^{2d})\,du = \operatorname{Beta}(1/(2d), 1 - 1/(2d))/(2d\pi).$$

If we can get a reasonable estimate $\hat{g}_\varepsilon$ of the $g_\varepsilon$, as in the method given below, then, by ignoring the term involving $\sigma^2$ in (3.1), we can arrive at the following criterion function

(3.2a) $$\phi(\nu, d) = SSE(\nu)/n + (1/n)\sum \hat{g}_\varepsilon(\pi j/n)/(1 + \nu s(\pi j/n)) \quad \text{or}$$

(3.2b) $$\phi(\nu, d) = SSE(\nu)/n + 2\nu^{-1/(2d)}\hat{g}_\varepsilon(0)c_6(d).$$

The first criterion given in (3.2a) is preferable, since the second one in (3.2b) is an approximation to the first when both $n$ and $\nu$ are large.

We will now concentrate on how to find a reasonable estimate of the spectral density of the error process $\{\varepsilon_t\}$ at zero. Let us assume that the error sequence is $AR(p)$, an autoregressive process of order $p$, where $p$ is unknown and needs to be estimated. Consider the local linear estimator $\hat{\mu}_t^+$ of $\mu_t$ on the basis of observations $Y_{t-k}, \ldots, Y_{t+k}$, where $k$ is approximately equal to $\sqrt{n}/2$. If $t < k+1$, then the estimator is based on $Y_1, \ldots, Y_{t+k}$, and a similar modification is done when $t > n - k$. Using a selection criterion such as AIC or BIC we can select an autoregressive model using the estimated error sequence $\widetilde{\varepsilon}_t = Y_t - \hat{\mu}_t^+$. Let $\hat{g}_\varepsilon(0)$ be the estimate of the spectral density of the error sequence at zero on the basis of the estimated autoregressive process.

A second approach is to minimize the innovations log likelihood for the state space model defined by the observation equation (1.1) and the state equations (1.2) and

$$\varepsilon_t = \phi_1 \varepsilon_{t-1} + \phi_2 \varepsilon_{t-2} + \delta_t, \tag{3.3}$$

where $\{\delta_t\}$ are mean zero i.i.d. with variance $\delta_d^2$. Again, AIC or BIC can be used to estimate the order $p$ of the autoregressive process. The estimators for $\mu$, say $\hat{\mu}$, are the usual Kalman smoothers, produced as a by-product when using the EM algorithm to estimate the unknown parameters $\phi_1$, $\phi_2$ and $\sigma_\delta^2$. The Kalman smoothers also produce the estimated mean square error of $E\|\hat{\mu} - \mu\|^2$, which can be used to set pointwise uncertainty limits for the smoothed trend.

We have used both methods above to obtain an estimate of the trend of the global temperature data, as given in Jones et al. (2000). Figure 1 shows the yearly average of land and marine temperature stations beginning in 1856 and ending in 2000. We have chosen a relatively long time span that indicates that the assumption of linearity, often made on the basis of temperature series beginning in 1900, may not be realistic over the long term. The first model selection criterion described above selected an AR(2) model for the error process $\{\varepsilon_t\}$ with parameter estimates $\hat{\phi}_1 = 0.3784$, $\hat{\phi} = -0.1660$, $\hat{\sigma}_\delta^2 = 0.0096$. The selected order of difference and the penalty parameters turned out be $\hat{d} = 2$ and $\hat{\nu} = 219.8$. Applying maximum likelihood (the second approach) yielded comparable values $\hat{\phi}_1 = 0.3882$, $\hat{\phi}_2 = -0.1641$ and $\hat{\sigma}_\delta^2 = 0.0095$. As a matter of fact, the trend estimates for these methods turned out to be indistinguishable. The fitted trend and data are plotted in Figure 1, and we note that the estimated trend conforms more to a nonlinear function with two periods of relative stable global temperatures and the two periods of rather steep increases, the last beginning at about 1975.

3.1. *Simulations.* We have done simulations in order to check the suitability of our criterion for sample sizes $n = 100$ and $n = 300$ with different signal-to-noise ratios (SNR). We have tried two cases when the true value of $d$ is either 1 or 2. So, the model we have tried is

$$Y_t = \mu_t + \varepsilon_t, \qquad t = 1, \ldots, n,$$

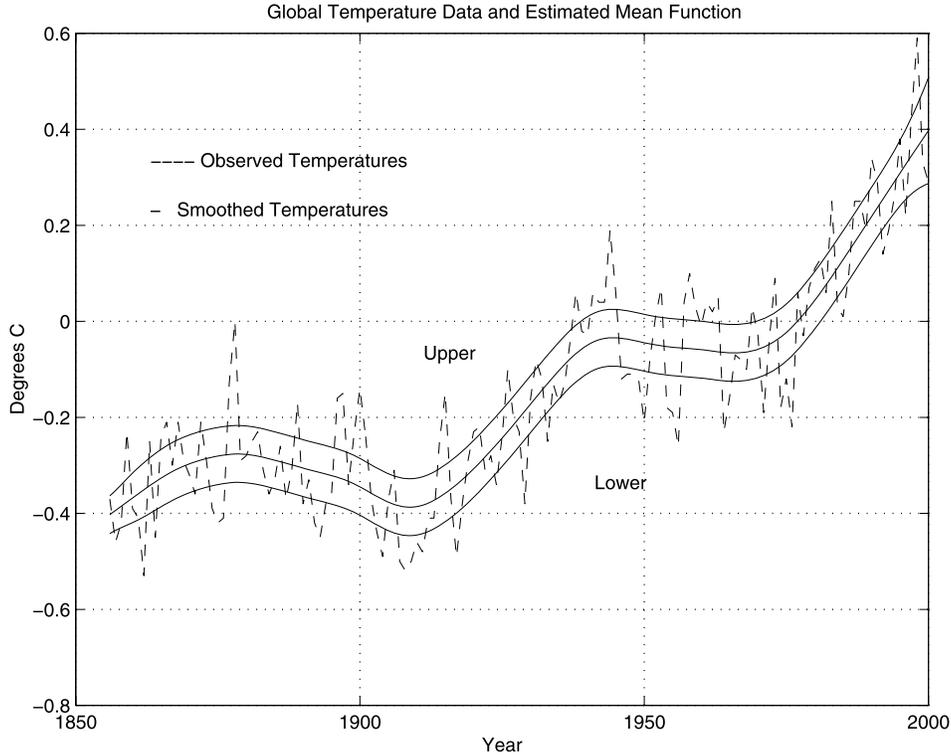

Fig. 1. *Yearly temperature anamalies (1856–2000) in degrees centigrade relative to the 1961–1990 mean. Solid lines are the fitted trend and the upper and lower point-wise 95% posterior probability points based on assuming normality of $\hat{\mu} - \mu$ and parameters fixed at their likelihood estimates.*

where $\{\varepsilon_t\}$ are i.i.d. $N(0,1)$, $\mu_t = \sum_{1 \leq j \leq t} (-1)^{t-j} \binom{-d}{t-j} \gamma_j$, where $\{\gamma_j\}$ are i.i.d. $N(0, \sigma_\gamma^2)$. We choose $\sigma_\gamma^2$ in such a way that the signal-to-noise ratio is either 2, 5 or 9. We have used the criterion given in (3.2a) and (3.2b) in order to estimate $\nu$ and $d$. Since the errors in the simulation are taken to be i.i.d., there are many ways to estimate its variance. We have used a fairly simple estimate here. The estimate used here is half times the average of the squares of the first differences of the observations. We have calculated the ratio $R = \inf_{d,\nu} \|\hat{\mu}(d, \nu) - \mu\| / \|\hat{\mu}(\hat{d}, \hat{\nu}) - \mu\|$ where the minimum in the numerator is over $\nu > 0$ and $d$ in $\{1, 2\}$, and $(\hat{\nu}, \hat{d})$ are obtained by minimizing the criterion function discussed above. We have calculated the mean, standard deviation and median of $R$ for various combinations of $d$ and SNR. All of the calculations are based on 400 repeats. How well we can estimate the underlying trend depends on the signal to noise ratio. The higher the SNR value, the better the estimate. The simulation results given in Table 1 show that the estimation methods proposed here work reasoanably well.

TABLE 1
*Simulated values of the performance ratio R*

|       |         | $n = 100$       |        | $n = 300$       |        |
|-------|---------|-----------------|--------|-----------------|--------|
|       |         | Mean(SD)        | Median | Mean(SD)        | Median |
| $d=1$ | SNR = 2 | 0.9381 (0.0603) | 0.9532 | 0.9599 (0.0387) | 0.9698 |
|       | SNR = 5 | 0.9411 (0.0488) | 0.9511 | 0.9726 (0.0306) | 0.9830 |
|       | SNR = 9 | 0.9367 (0.0460) | 0.9435 | 0.9728 (0.0285) | 0.9830 |
| $d=2$ | SNR = 2 | 0.7833 (0.1839) | 0.8231 | 0.8343 (0.1582) | 0.8768 |
|       | SNR = 5 | 0.8175 (0.1531) | 0.8502 | 0.8424 (0.1548) | 0.8923 |
|       | SNR = 9 | 0.8377 (0.1451) | 0.8735 | 0.8696 (0.1335) | 0.9142 |

## APPENDIX

We will begin this section with a number of notations and definitions which will be used in the proofs.

NOTATION 1. For any square matrix $A$ of order $n$, we will denote its singular values by $\sigma_1(A), \ldots, \sigma_n(A)$ and its eigenvalues by $\lambda_1(A), \ldots, \lambda_n(A)$ (singular values of $A$ are the positive square roots of the eigenvalues of $A'A$).

NOTATION 2. The function $1 - e^{iu}$, $-\pi \leq u \leq \pi$, will be dented by $s_0(u)$ and $|s_0(u)|^2 = 2 - 2\cos(u)$ will be denoted by $s(u)$.

It is important to note that the $d$th order finite difference matrix as given in (2.3) is a Toeplitz matrix with symbol $s_0^d$. The Toeplitz matrices asociated with $s_0^d$ and $s^d$ will be used often in the proofs. There are a number of different matrix norms that will come into play, and we will define them here.

DEFINITION 1. Let $\sigma_j(A)$, $j = 1, \ldots, n$, be the singular values of the matrix $A$. The following three norms are used widely:

(i) Spectral radius norm. $\|A\| = \max \sigma_j(A)$;
(ii) Trace norm. $\|A\|_1 = \sum \sigma_j(A)$;
(iii) Frobenius norm. $\|A\|_2 = \{\sum \sigma_j^2(A)\}^{1/2}$.

DEFINITION 2. A matrix $T = ((b_{jk}))$ is of Toeplitz type if $b_{jk} = b_{j-k}$. If $b_{j-k}$ is given by $\int_{-\pi}^{\pi} \exp(i(j-k)u) f(u)\, du/(2\pi)$, then the function $f$ is called the symbol of the Toeplitz matrix and often the notation $T(f)$ is used to denote the Toeplitz matrix. The submatrix consisting of the first $n$ rows and columns will be denoted by $T_n(f)$.

DEFINITION 3. A matrix $H = ((b_{jk}))$ is of Hankel type if $b_{jk}$ is of the form $b_{j+k}$. If $b_{j+k}$ is given by $\int_{-\pi}^{\pi} \exp(i(j+k)u) f(u) \, du/(2\pi)$, then $f$ is called the symbol of the Hankel matrix and often the notation $H(f)$ is used to denote the Hankel matrix. The submatrix consisting of the first $n$ rows and columns will be denoted by $H_n(f)$.

We will now define another special type of matrices called the circulants. The circulants can be used to approximate Toeplitz matrices. Let $P_0$ be the $n \times n$ cyclical permutation matrix whose element $(j,k)$ is 1 if $j-k = 1[\mathrm{mod}\, n]$ and 0 otherwise. Then, $P_0$ has the form given below:

$$
(A.1) \qquad P_0 = \begin{pmatrix} 0 & 0 & 0 & \cdot & \cdot & \cdot & \cdot & 1 \\ 1 & 0 & 0 & \cdot & \cdot & \cdot & \cdot & 0 \\ 0 & 1 & 0 & \cdot & \cdot & \cdot & \cdot & 0 \\ \cdot & \cdot & \cdot & \cdot & \cdot & \cdot & \cdot & \cdot \\ \cdot & \cdot & \cdot & \cdot & \cdot & \cdot & \cdot & \cdot \\ \cdot & \cdot & \cdot & \cdot & \cdot & \cdot & \cdot & \cdot \\ \cdot & \cdot & \cdot & \cdot & \cdot & \cdot & 0 & 0 \\ 0 & \cdot & 0 & \cdot & \cdot & \cdot & 1 & 0 \end{pmatrix}.
$$

DEFINITION 4. A square matrix $C_n = ((b_{jk}))$ of order $n$ is called a circulant if $b_{jk} = b_{j-k[\mathrm{mod}\, n]}$. If $b_j$'s are given by $f(u) = -\sum_{-r \le j \le r} b_j e^{iju}$, then the circulant $C_n$ is said to be generated by the symbol $f$, and we write $C_n(f)$ to denote it. If $P_0$ is the cyclical permutation matrix as given above in (A.1), then we can express $C_n(f) = \sum_{-r \le j \le r} b_j P_0^j$.

We will first write down a few important lemmas, the proofs of which will be given after those of Theorems 1 and 2. We begin with an interlacing theorem due to Weyl [see Theorem 4.3.6 in Horn and Johnson (1985)].

THEOREM 6. *Let $A$ and $B$ be two real symmetric matrices so that $A - B$ has rank at most $r$. Let $\lambda_1(A) \le \cdots \le \lambda_n(A)$ and $\lambda_1(B) \le \cdots \le \lambda_n(B)$ be the eigenvalues of the matrices $A$ and $B$. Then:*

(a) $\lambda_j(A) \le \lambda_{j+r}(B)$, $j = 1, \ldots, n-r$,
(b) $\lambda_j(B) \le \lambda_{j+r}(A)$, $j = 1, \ldots, n-r$.

The next result on singular values, which is analogous to the previous one, follows from the result given on page 423 in Horn and Johnson (1985).

THEOREM 7. *Let $A$ and $B$ be two matrices of order $n \times n$, and the matrix $A - B$ has rank at most $r$. Let $\sigma_1(A) \ge \sigma_2(A) \ge \cdots$ be the singular values of $A$, and, similarly, let $\sigma_1(B) \ge \sigma_2(B) \ge \cdots$ be the singular values of $B$. Then:*

(a) $\sigma_{j+r}(A) \leq \sigma_j(B)$, $j = 1, \ldots, n - r$,
(b) $\sigma_{j+r}(B) \leq \sigma_j(A)$, $j = 1, \ldots, n - r$.

The next lemma finds the bounds on the singular values and the trace norm of a Hankel matrix.

LEMMA 1. *Let $H = ((b_{j+k}))$ be a $n \times n$ real Hankel matrix. Let $\sigma_1(H) \geq \cdots \geq \sigma_n(H)$ be the singular values of $H$. Then:*

(a) $\sigma_j(H) \leq \sum_{j+1 \leq t \leq 2n} |b_t|$,
(b) $\|H\|_1 = \sum_{1 \leq j \leq n} |\sigma_j(H)| \leq 2 \sum_{1 \leq l \leq 2n} l|b_l|$.

The following important lemma tells us how a Toeplitz matrix can be approximated by a circulant matrix.

LEMMA 2. *Let $T_n(f)$ be a Toeplitz matrix with the symbol $f(u) = \sum_{-N \leq j \leq N} b_j e^{-iju}$, where $N < n/2$. Let $C_n(f)$ be the circulant matrix given by $\sum_{-N \leq j \leq N} b_j P_0^j$, where the generator permutation matrix is as given in (A.1). For $1 \leq j \leq n$, consider the vector $e_j$ whose $t$th element is $n^{-1/2} \exp(-i2\pi jt/n)$, $1 \leq t \leq n$. Then, $e_1, \ldots, e_n$ are orthonormal. The following results hold:*

(a) *When $b_j = b_{-j}$, the eigenvalues (unordered) and the corresponding eigenvectors of the circulant $C_n(f)$ are given by $f(2\pi j/n)$ with the corresponding eigenvectors $e_j$;*
(b) $\operatorname{rank}(C_n(f) - T_n(f)) \leq 2N$;
(c) $\|C_n(f) - T_n(f)\|_1 \leq 2 \sum_{1 \leq j \leq N} (j+1)|b_j|$;
(d) *The circulant matrix $C_n(f)$ can written as $\sum_{1 \leq j \leq n} f(2\pi j/n) e_j e_j^*$.*

The next result compares the sum of squares of the singular values of a matrix $A$ to the sum of the singular values of a another matrix $B$ when the matrix $A - B$ is of finite rank. The proof is omitted as it an easy consequence of the interlacing theorem.

LEMMA 3. *Let $A$ and $B$ be two square matrices of order $n$, and the rank of matrix $A - B$ does not depend on $n$. Then,*

$$\sum_{1 \leq j \leq n} \sigma_j^2(A) - \sum_{1 \leq j \leq n} \sigma_j^2(B) = O(1)(\sigma_1^2(A) + \sigma_1^2(B)).$$

We will present a few known results without proofs. These results will be useful in our proofs. The first result [Theorem 1.1 in Böttcher and Grudsky (2000)] obtains an upper bound of the norm of a Teeplitz matrix in terms of the supremum norm of its symbol.

THEOREM 8. *Let $T_n(f)$ be a Toeplitz matrix with symbol $f$. Let $|f|_\infty$ be the supremum norm of $f$. Then,*
$$\|T_n(f)\| \leq |f|_\infty.$$

The following result tells us that the matrix $\overline{S}'_{-d}\overline{S}_{-d}$ does not differ much from the Toeplitz matrix $T_n(s^d)$.

LEMMA 4. *All the elements of the matrix $\overline{S}'_{-d}\overline{S}_{-d} - T_n(s^d)$ are zero except for the first and the last principal submatrices of order $d$.*

PROOF OF THEOREM 1. In the proofs, we will assume that, for any real symmetric matrix $C$ of order $n$, its eigenvalues denoted by $\lambda_1(C), \ldots, \lambda_n(C)$ are ordered from the smallest to the largest. Also note that the unordered eigenvalues of a circulant $C_n(f)$ are given by $f(2\pi j/n)$, $j = 1, \ldots, n$. Hence, it possible that $\lambda_j(C_n(f)) \neq f(2\pi j/n)$ for some (or even all) values of $j$. For notational simplicity, we will denote the matrix $\overline{S}'_{-d}\overline{S}_{-d}$ by $U$.

First, note that the estimate $\hat{\mu}$ is given by $(I + \nu \overline{S}'_{-d}\overline{S}_{-d})^{-1} Y = (I + \nu U)^{-1} Y$ [see (2.5)] and $\overline{\mu} = (I + \nu \overline{S}'_{-d}\overline{S}_{-d})^{-1}\mu = (I + \nu U)^{-1}\mu$. Since $g_\varepsilon$ is the spectral density function of the process $\{\varepsilon_t\}$, the variance–covariance matrix of $\{\varepsilon_t : t = 1, \ldots, n\}$ is given by the Toeplitz matrix $T_n(g_\varepsilon)$. Hence, we have

(A.2) $$E\|\hat{\mu} - \overline{\mu}\|^2/n = \operatorname{tr}((I + \nu U)^{-2} T_n(g))/n.$$

The main idea behind the proof is to use approximate $\overline{S}'_{-d}\overline{S}_{-d} = U$ and $T_n(g_\varepsilon)$ by circulants and then use the well-known theory of circulants to get the result.

Recall that the spectral density function of the error process $\{\varepsilon_t\}$ is given by $g_\varepsilon(u) = \sum_{-\infty < t < \infty} \rho_\varepsilon(t) e^{-itu}$, where the covariances satisfy the condition $\sum_{1 \leq t < \infty} t|\rho_\varepsilon(t)| < \infty$.

Let $N = [n/4]$, the integer part of $n/4$, and define $\overline{g}_{\varepsilon N}(u) = \sum_{|t| \leq N} \rho_\varepsilon(t) e^{-itu}$. Then,
$$\|\overline{g}_{\varepsilon N} - g_\varepsilon\|_\infty \leq \sum_{|t| > N} |\rho_\varepsilon(t)| \leq N^{-1} \sum_{|t| > N} t|\rho_\varepsilon(t)| = o(1/n).$$

Now, define $g_{\varepsilon N}(u) = \max(\overline{g}_{\varepsilon N}(u), 0)$. Since $g_\varepsilon$ is a nonnegative function, $\|g_{\varepsilon N} - g_\varepsilon\|_\infty = o(1/n)$. From Theorem 8, we have
$$\|T_n(g_\varepsilon - g_{\varepsilon N})\| \leq \|g_\varepsilon - g_{\varepsilon N}\|_\infty = o(1/n).$$

Hence,
$$|\operatorname{tr}((I + \nu U)^{-2} T_n(g_\varepsilon)) - \operatorname{tr}((I + \nu U)^{-2} T_n(g_{\varepsilon N}))|$$
(A.3) $$= |\operatorname{tr}((I + \nu U)^{-2} T_n(g_\varepsilon - g_{\varepsilon N}))|$$
$$\leq \operatorname{tr}((I + \nu U)^{-2})\|T_n(g_\varepsilon - g_{\varepsilon N})\| = o(1/n)\operatorname{tr}((I + \nu U)^{-2}) = o(1).$$

Now, using part (c) of Lemma 2, we get

(A.4)
$$\begin{aligned}|\text{tr}((I+\nu U)^{-2}C_n(g_{\varepsilon N})) &- \text{tr}((I+\nu U)^{-2}T_n(g_{\varepsilon N}))| \\ &\leq \|(I+\nu U)^{-2}\|\|C_n(g_{\varepsilon N})-T_n(g_{\varepsilon N})\|_1 \\ &\leq O(1)\sum_{1\leq t\leq N}(t+1)|\rho_\varepsilon(t)| < \infty.\end{aligned}$$

So, combining (A.3) and (A.4), we get

(A.5) $$\text{tr}((I+\nu U)^{-2}T_n(g_\varepsilon)) = \text{tr}((I+\nu U)^{-2}C_n(g_{\varepsilon N})) + O(1).$$

Note that $g_\varepsilon$ is differentiable because of the assumption $\sum_{1\leq t<\infty} t|\rho_\varepsilon(t)| < \infty$. Since $g_\varepsilon$ is assumed to be positive on $[\pi, \pi]$, we can assume that $g_{\varepsilon N}$ is also positive as $N \to \infty$. Now, the eigenvalues (unordered) of the circulant $C_n(g_{\varepsilon N})$ are given by $g_{\varepsilon N}(2\pi j/n)$, and they are all nonnegative. Consequently, we can define a positive square root of the matrix $C_n(g_{\varepsilon N})$, and it is given by $C_n(g_{\varepsilon N})^{1/2} = \sum g_{\varepsilon N}(2\pi j/n)^{1/2} e_j e_j^*$, where $e_j$'s are the eigenvectors given in Lemma 2.

Consider the matrices

(A.6)
$$\begin{aligned}B &= (I+\nu U)^{-1}C_n(g_{\varepsilon N})^{1/2} \quad \text{and} \\ A &= (I+\nu C_n(s^d))^{-1}C_n(g_{\varepsilon N})^{1/2},\end{aligned}$$

where $s(u) = 2 - 2\cos u$ is as defined in Notation 2 at the beginning of this section.

So, from (A.5) and (A.6), we have

(A.7) $$\text{tr}((I+\nu U)^{-2}T_n(g_\varepsilon)) = \text{tr}(B'B) + O(1) = \sum_{1\leq j\leq n}\sigma_j(B)^2 + O(1).$$

Now, note that

$$A = B + \nu(I+\nu C_n(s))^{-1}(U - C_n(s))(I+\nu C_n(g_{\varepsilon N}))^{-1}C_n(g_{\varepsilon N})^{1/2}.$$

Lemma 4 tells us $T_n(s^d) - U$ has rank at most $2d$. Since $T_n(s^d)$ is a banded matrix, by part (b) of Lemma 2 we see that the rank of $C_n(s^d) - T_n(s^d)$ is at most $2d$. So, the rank of $C_n(s^d) - U$ is at most $4d$. Consequently, the rank of the matrix $A - B$ is no larger than $4d$. Now, note that eigenvalues of the matrix $A'A = C_n(g_{\varepsilon N})^{1/2}(I+\nu C_n(s^d))^{-2}C_n(g_{\varepsilon N})^{1/2}$ and $B'B = C_n(g_{\varepsilon N})^{1/2}(I+\nu U)^{-2}C_n(g_{\varepsilon N})^{1/2}$ are bounded above by a positive constant, which is independent of $n$ and $\nu$. Hence, by Lemma 3,

$$\text{tr}(B'B) = \sum_{1\leq j\leq n-4d}\sigma_j(A)^2 + O(1).$$

Note that the unordered eigenvalues of the matrix $A'A$ are $(1+\nu s(2\pi j/n)^d)^{-2} \times g_{\varepsilon N}(2\pi j/n)$, $j = 1,\ldots,n$. Hence, we can conclude that

$$\text{tr}(B'B) = \sum_{1 \le j \le n} \sigma_j(A)^2 + O(1)$$
(A.8)
$$= \sum_{1 \le j \le n} g_{\varepsilon N}(2\pi j/n)(1+\nu s(2\pi j/n)^d)^{-2} + O(1).$$

Note that $g_\varepsilon$, and, hence, $g_{\varepsilon N} I(g_{\varepsilon N} > 0)$ are differentiable with bounded first derivatives because of the assumption $\sum_{1 \le t < \infty} t|\rho_\varepsilon(t)| < \infty$, and that $\|g_\varepsilon - g_{\varepsilon N}\|_\infty = o(n^{-1})$. From (A.8), we get

$$\text{tr}(B'B) = \sum_{1 \le j \le n} g_{\varepsilon N}(2\pi j/n)(1+\nu s(2\pi j/n)^d)^{-2} + O(1)$$

$$= \int_0^n g_\varepsilon(2\pi u/n)(1+\nu s(2\pi u/n)^d)^{-2}\, du + O(1)$$

$$= n(2\pi)^{-1} \int_0^{2\pi} g_\varepsilon(u)(1+\nu s(u)^d)^{-2} + O(1).$$

From the last expression and (A.2) and (A.8), we get

$$E\|\hat{\mu} - \overline{\mu}\|^2/n = n^{-1}\text{tr}(B'B) + O(1/n)$$
(A.9)
$$= (2\pi)^{-1} \int_0^{2\pi} g(u)(1+\nu s(u)^d)^{-2} + O(1/n)$$

$$= \pi^{-1} \int_0^{\pi} g(u)(1+\nu s(u)^d)^{-2} + O(1/n).$$

Now, if we denote $\phi(u) = \{\sin(u/2)/(u/2)\}^{2d}$, then we can write

(A.10) $\quad s(u)^d = (2 - 2\cos u)^d = 4^d \sin(u/2)^{2d} = u^{2d}\phi(u).$

Making a change of variable $z = \nu^{1/(2d)} u$, we get

$$E\|\hat{\mu} - \overline{\mu}\|^2/n$$
(A.11)
$$= \pi^{-1}\nu^{-1/(2d)} \int_0^{\pi\nu^{1/(2d)}} g_\varepsilon(z\nu^{-1/(2d)})(1+z^{2d}\phi(z\nu^{-1/(2d)}))^{-2}\, dz$$
$$+ O(1/n).$$

Now, note $\phi(0) = 1$ and that the function $g_\varepsilon$ and $\phi$ have bounded derivatives on $[0, \pi]$. Consequently, the last expression yields

$$E\|\hat{\mu} - \overline{\mu}\|^2/n = (2\pi)^{-1}\nu^{-1/(2d)} g_\varepsilon(0) \int_0^\infty (1+z^{2d})^{-2}\, dz(1 + O(\nu^{-1/(2d)}))$$
$$+ O(1/n).$$

The proof of this result is now complete, once we note that

$$\int_0^\infty (1+z^{2d})^{-2}\,dz = \text{Beta}(1/(2d), 2 - 1/(2d))/(2d). \qquad \square$$

PROOF OF THEOREM 2. As in the proof of Theorem 1, here, too, we will assume that, for any real symmetric matrix $C$ of order $n$, its eigenvalues, denoted by $\lambda_1(C), \ldots, \lambda_n(C)$, are ordered from the smallest to the largest. Unordered eigenvalues of a circulant $C_n(f)$ are given by $f(2\pi j/n)$, $j = 1, \ldots, n$, and it is possible that $\lambda_j(C_n(f)) \neq f(2\pi j/n)$ for some (or even all) values of $j$. As in the proof of the last theorem, we will denote $\overline{S}'_{-d}\overline{S}_{-d}$ by $U$.

A matrix representation of $\hat{\mu}$ is given in (2.5). As a consequence, we have

$$\overline{\mu} - \mu = (I + \nu U)^{-1}\mu - \mu = -\nu(I + \nu U)^{-1}U\mu,$$

since

$$\overline{S}_{-d}\mu = (\nabla^d \mu_{d+1}, \ldots, \nabla^d \mu_n)' = (\gamma_{d+1}, \ldots, \gamma_n)',$$

where $\gamma_j$ are i.i.d. with mean zero and variance $\sigma_\gamma^2$. Hence,

$$(A.12) \qquad E\|\overline{\mu} - \mu\|^2 = \nu^2 \sigma_\gamma^2 \operatorname{tr}((I + \nu U)^{-2}U).$$

Let $\psi(u) = u/(1 + \nu u^2)$, when $u$ is a real number. Now, we can write the matrix $(I + \nu U)^{-2}U$ as $\psi(U)$. Hence, a re-expression of the relation (A.12) is given by

$$(A.13) \qquad E\|\overline{\mu} - \mu\|^2 = \nu^2 \sigma_\gamma^2 \operatorname{tr}(\psi(U)) = \nu^2 \sigma_\gamma^2 \sum_{1 \le j \le n} \psi(\lambda_j(U)),$$

where $\lambda_1(U), \ldots, \lambda_n(U)$ are the eigenvalues of the matrix $U$. As in the proof of Theorem 1, we can now approximate $U$ by the circulant matrix $C_n(s^d)$. Now, note that, for any $j = 1, \ldots, n$, both $\psi(\lambda_j(U))$ and $\psi(\lambda_j(C_n(s^d)))$ are bounded above by $\nu^{-1}$. Since $\psi$ is an increasing function on $[0, \infty)$, we can follow an argument similar to the one given in the proof of Theorem 1 to show that

$$\sum_{1 \le j \le n} \psi(\lambda_j(U)) = \sum_{1 \le j \le n} \psi(\lambda_j(C_n(s^d))) + O(\nu^{-1}).$$

Since the unordered eigenvalues of $C_n(s^d)$ are $s(2\pi j/n)^d$, $j = 1, \ldots, n$, by Lemma 2, from the last expression we have

$$(A.14) \qquad \sum_{1 \le j \le n} \psi(\lambda_j(U)) = \sum_{1 \le j \le n} \psi(s(2\pi j/n)^d) + O(\nu^{-1}).$$

Noting that the function $\psi$ is bounded on $[0,\infty)$, we have

$$
\begin{aligned}
\sum_{1\le j\le n} \psi(s(2\pi j/n)^d) &= \int_0^n \psi(s(2\pi u/n)^d)\,du + O(\nu^{-1}) \\
&= n(2\pi)^{-1}\int_0^{2\pi} \psi(s(u)^d)\,du + O(\nu^{-1}) \\
&= n\pi^{-1}\int_0^{\pi} \psi(s(u)^d)\,du + O(\nu^{-1}).
\end{aligned}
\tag{A.15}
$$

From the relations (A.13), (A.14) and (A.15), we get

$$
E\|\overline{\mu}-\mu\|^2/n = \nu^2\sigma_\gamma^2\pi^{-1}\int_0^{\pi}\psi(s(u)^d)\,du + O(\nu/n^{2d}). \tag{A.16}
$$

If we write $s(u)^d = u^{2d}\phi(u)$ as in (A.10), then, by a change of variable $z = \nu^{1/(2d)}u$, we get

$$
\begin{aligned}
\int_0^{\pi}\psi(s(u)^d)\,du &= \int_0^{\pi} u^{2d}\phi(u)(1+\nu u^{2d}\phi(u))^{-2}\,du \\
&= \nu^{-1-1/(2d)}\int_0^{\pi\nu^{1/(2d)}} z^{2d}\phi(z\nu^{-1/(2d)}) \\
&\qquad\times(1+z^{2d}\phi(z\nu^{-1/(2d)}))^{-2}\,dz.
\end{aligned}
\tag{A.17}
$$

Since $\phi(0)=1$ and $\phi$ has a bounded first derivative, calculations will show that

$$
\begin{aligned}
\int_0^{\pi\nu^{1/(2d)}} &z^{2d}\phi(z\nu^{-1/(2d)})(1+z^{2d}\phi(z\nu^{-1/(2d)}))^{-2}\,dz \\
&= \int_0^{\infty} z^{2d}(1+z^{2d})^{-2}\,dz\,(1+O(\nu^{-1/(2d)})).
\end{aligned}
\tag{A.18}
$$

Combining (A.16), (A.17) and (A.18), we get

$$
\begin{aligned}
E\|\overline{\mu}-\mu\|^2/n &= \nu^{1-1/(2d)}\sigma_\gamma^2\pi^{-1}\int_0^{\infty} z^{2d}(1+z^{2d})^{-2}\,dz\,(1+O(\nu^{-1/(2d)})) \\
&\quad + O(\nu/n^{2d}).
\end{aligned}
$$

This completes the proof of this result, once we note that

$$
\int_0^{\infty} z^{2d}(1+z^{2d})^{-2}\,dz = \mathrm{Beta}(1+1/(2d),1-1/(2d))/(2d). \qquad \square
$$

PROOF OF THEOREM 3. From the proof of Theorem 2, we have

$$
\overline{\mu}-\mu = -\nu(I+\nu U)^{-1}U\mu = -\nu(I+\nu U)^{-1}U(\mu_{1t}+S_{d_0}\gamma), \tag{A.19}
$$

where $U = \overline{S}'_{-d}\overline{S}_{-d}$.

(a) We will assume that the polynomial part is written in the form $\mu_{1t} = \sum_{1 \leq j \leq d_0 - 1} \beta_j (t/n)^j$ and that the coefficients $\beta_j$ are constants. The proof for the case when $\beta_j$ are random with finite mean and variance is the same. Note that

$$
\begin{aligned}
E\{\|\overline{\mu} - \mu\|^2\}/n &= n^{-1}\nu^2[\|(I + \nu U)^{-1} U \mu_1\|^2 \\
&\quad + \sigma_\gamma^2 \operatorname{tr}((I + \nu U)^{-1} U S_{d_0} S'_{d_0} U (I + \nu U)^{-1})] \\
&= n^{-1}\nu^2[\|(I + \nu U)^{-1} U \mu_1\|^2 + \tau^2 n^{-2d_0} \operatorname{tr}(AA')],
\end{aligned}
\tag{A.20}
$$

where $A = (I + \nu U)^{-1} U S_{d_0}$. Note that $\|\overline{S}_{-d}\mu_1\| = O(n^{-d})$ and, hence,

$$
\begin{aligned}
\|(I + \nu U)^{-1} U \mu_1\|^2 &= O(n^{-2d}) \|(I + \nu \overline{S}'_{-d}\overline{S}_{-d})^{-1} \overline{S}'_{-d}\| \\
&= O(n^{-2d}\nu^{-1}).
\end{aligned}
\tag{A.21}
$$

Since $\overline{S}_{-d}S_d$ is a $(n - d) \times n$ is a parttioned matrix of the form $[0 : I] = \widetilde{I}$, where the first matrix in the partition is a $(n-d) \times d$ matrix of zeros and the second matrix is the identity matrix of order $n - d$. So the matrix $\overline{S}_{-d}S_{d_0}$ can be rewritten as $\overline{S}_{-d}S_d S_{d_0 - d} = \widetilde{I} S_{d_0 - d}$. It is known that the largest singular value of $S_{d_0 - d}$ is of order $n^{d_0 - d}$ [see Theorem 2 in Burman (2006)]. Hence,

$$
\begin{aligned}
\operatorname{tr}(AA') &\leq \operatorname{tr}((I + \nu \overline{S}'_{-d}\overline{S}_{-d})^{-1} \overline{S}'_{-d}\overline{S}_{-d}(I + \nu \overline{S}'_{-d}\overline{S}_{-d})^{-1}) \|\widetilde{I} S_{d_0 - d}\|^2 \\
&= \operatorname{tr}((I + \nu \overline{S}'_{-d}\overline{S}_{-d})^{-1} \overline{S}'_{-d}\overline{S}_{-d}(I + \nu \overline{S}'_{-d}\overline{S}_{-d})^{-1}) O(n^{2(d_0 - d)}).
\end{aligned}
$$

The proof of Theorem 2 [(A.14) through (A.2.8)] shows that

$$
\operatorname{tr}((I + \nu \overline{S}'_{-d}\overline{S}_{-d})^{-1} \overline{S}'_{-d}\overline{S}_{-d}(I + \nu \overline{S}'_{-d}\overline{S}_{-d})^{-1}) = O(n\nu^{-1 - 1/(2d)}).
$$

Hence, we get $\operatorname{tr}(AA') = O(n^{-2d}\nu^{-1 - 1/(2d)})$. Now, combining this result with those from (A.20) and (A.21), we get

$$
E\{\|\overline{\mu} - \mu\|^2\}/n = O(n^{-2d}\nu^{-1}) + O(n\nu^{-1 - 1/(2d)}) = O(n^{-2d}\nu^{-1}).
$$

(b) Since $\overline{S}_{-d}\mu_1 = 0$, from (A.20) we get

$$
E\{\|\overline{\mu} - \mu\|^2\}/n = n^{-1}\nu^2 \tau^2 n^{-2d_0} \operatorname{tr}(AA').
$$

Let

$$
C_n(s_0^{d_0})^- = \sum_{1 \leq j \leq n - 1} s_0(2\pi j/n)^{-d_0} e_j e_j^*.
$$

Then, $C_n(s_0^{d_0})^-$ is generalized inverse of $C_n(s_0^{d_0})$. We will approximate $AA'$ by $BB'$, where $B = (I + \nu C_n(s^d))^{-1} C_n(s^d) C_n(s^{d_0/2})^-$. Since the rank of $S'_{-d_0}S_{-d_0} - C_n(s^{d_0})$ is no larger than $2d_0$, and the rank of $C_n(s_0^{d_0})$ is $n - 1$, the rank of $S_{d_0} - C_n(s_0^{d_0})^- = S_{-d_0}^{-1} - C_n(s_0^{d_0})^-$ is at most $2d_0 + 1$. Also, note that the rank of $U - C_n(s^d)$ is at most $2d + 1$. So, using Lemma 3, we get

$$E\{\|\bar{\mu} - \mu\|^2\}/n = \tau^2 \nu^2 n^{-2d_0} \operatorname{tr}(AA')$$
$$= \tau^2 \nu^2 n^{-2d_0} \operatorname{tr}(BB^*) + O(\nu^2 n^{-2d_0})(\sigma_1(A)^2 + \sigma_1(B)^2).$$

We will later show that, in the last expression, the second term involving $\sigma_1(A)^2$ and $\sigma_1(B)^2$ is small in comparison to the first term; that is,

$$E\{\|\bar{\mu} - \mu\|^2\}/n = \tau^2 \nu^2 n^{-2d_0} \operatorname{tr}(BB^*)[1 + o(1)].$$

Note that the smallest eigenvalue of $BB^*$ is zero and the rest of the eigenvalues (unordered) are given by $s(2\pi j/n)^{2d-d_0}/(1+\nu s(2\pi j/n)^d)^2$, $j = 1, \ldots, n-1$.

So,

$$\tau^2 \nu^2 n^{-2d_0} \operatorname{tr}(BB^*) = \tau^2 n^{-2d_0} \sum_{1 \le j \le n-1} s(2\pi j/n)^{2d-d_0}/(1+\nu s(2\pi j/n)^d)^2.$$

When $d \ge d_0$, an argument similar to the one used in the proof of Theorem 2 will show that

$$\tau^2 \nu^2 n^{-2d_0} \operatorname{tr}(BB^*)$$
$$= \tau^2 (\nu^{1/(2d)}/n)^{2d_0 - 1}$$
$$\quad \times \operatorname{Beta}((2d_0 - 1)/(2d), 2 - (2d_0 - 1)/(2d))/(2d\pi)[1 + o(1)].$$

What is left to show is that

$$\nu^2 n^{-2d_0}(\sigma_1(A)^2 + \sigma_1(B)^2) = o(1)(\nu^{1/(2d)}/n)^{2d_0 - 1}.$$

We will first prove the case for $\sigma_1(B)^2$. Recall that the smallest eigenvalue of $BB^*$ is zero, and the rest of the eigenvalues (unordered) are given by $\psi_j = s(2\pi j/n)^{2d-d_0}/(1+\nu s(2\pi j/n)^d)^2$, $j = 1, \ldots, n-1$. Since the largest eigenvalue of $BB^*$ is no larger that $\nu^{-2+d_0/d}$, we have that

$$\nu^2 n^{-2d_0} \sigma_1(B)^2 \le (\nu^{1/(2d)})^{2d_0} = o(1)(\nu^{1/(2d)}/n)^{2d_0 - 1}.$$

Now, let us find the bound for the term $\sigma_1(A)^2$. Let $F = S'_{-1} S_{-1}$. Calculations will show that

$$\overline{S}_{-d} S_{d_2} S'_{d_0} \overline{S}'_{-d} = \overline{S}_{-d}(S_1 S'_1)^{d_0} \overline{S}'_{-d} = \overline{S}_{-d} F^{-d_0} \overline{S}'_{-d}.$$

It can be shown that $U = \overline{S}'_{-d} \overline{S}_{-d} \le (S'_{-1} S_{-1})^d = F^d$, where, for any two matrices, the notation "$C \le D$" means that $D - C$ is nonnegative definite. If $C$ and $D$ are nonnegative definite and $C \le D$, then it can be shown that $(I+C)^{-1}C \le (I+D)^{-1}D$. Consequently, the largest eigenvalue of $AA'$ is no larger than the largest eigenvalue of $(I + \nu F^d)^{-1} F^d F^{-d_0} F^d (I + \nu F^d) = (I + \nu F^d)^{-2} F^{2d-d_0}$. By Gershgorin's result, one can see that all the eigenvalues of $F$ are bounded above by 4. So, the largest eigenvalue of $(I + \nu F^d)^{-2} F^{2d-d_0}$ is bounded above by $\nu^{-2+d_0/d}$. Consequently, the largest eigenvalues of $AA'$ is no larger that $\nu^{-2+d_0/d}$. Hence, we conclude that

$$\nu^2 n^{-2d_0} \sigma_1(A)^2 \leq (\nu^{1/(2d)})^{2d_0} = o(1)(\nu^{1/(2d)}/n)^{2d_0-1}. \qquad \square$$

PROOF OF LEMMA 1. (a) Since the singular values of $H$ are the positive square root of the eigenvalues of $H'H$, by the Courant–Fischer minimax theorem [Theorem 7.3.10 in Horn and Johnson (1985)], the square of the $j$th singular value of $H$ is

$$\sigma_j^2(H) = \min_{u_1,\ldots,u_{j-1} \text{ in } R^n} \max_{x \in R^n, \|x\|=1, x \perp \{u_1,\ldots,u_{j-1}\}} x'H'Hx.$$

Now, take any vector $x$ in $R^n$ whose first $j-1$ components are zero. In the case we will be basically concerned with, the principal submatrix of $H'H$ consists of the last $n-j+1$ columns and rows. Applying Gersgorin's theorem on the localization of eigenvalues [see Theorem 6.1.1 in Horn and Johnson (1985)], the largest eigenvalue of this principal submatrix is no larger than

$$\max_{j \leq s \leq n} \sum_{j \leq t \leq n} \left| \sum_{1 \leq l \leq n} b_{s+l} b_{l+t} \right|.$$

So, we have

$$\sigma_j^2(H) \leq \max_{j \leq s \leq n} \sum_{j \leq t \leq n} \left| \sum_{1 \leq l \leq n} b_{s+l} b_{l+t} \right|$$

$$\leq \max_{j \leq s \leq n} \sum_{1 \leq l \leq n} |b_{s+l}| \left\{ \sum_{l+j \leq t \leq n+l} |b_t| \right\}$$

$$\leq \max_{j \leq s \leq n} \sum_{1 \leq l \leq n} |b_{s+l}| \left\{ \sum_{l+j \leq t \leq n+l} |b_t| \right\}$$

$$\leq \max_{j \leq s \leq n} \sum_{s+1 \leq l \leq n+s} b_l \left( \sum_{j+1 \leq t \leq 2n} |b_t| \right) \leq \left( \sum_{j+1 \leq t \leq 2n} |b_t| \right)^2.$$

(b) Using part (a), we have

$$\sum_{1 \leq j \leq n} |\sigma_j(H)| \leq \sum_{1 \leq j \leq n} \sum_{j+1 \leq t \leq 2n} |b_t| \leq \sum_{1 \leq t \leq 2n} t|b_t|$$

and this completes the proof. $\square$

PROOF OF LEMMA 2. (a) This part follows from the well-known results on circulant matrices [see Chapter 4 in Marcus and Minc (1992) or Tyrtyshnikov (1996)].

(b) and (c). Let $W_n$ be the orthogonal matrix that has the property that, for any $x = (x_1, \ldots, x_n)'$ in $R^n$, it flips it indexwise; that is, $W_n x = (x_n, \ldots, x_1)'$. In turns out that $W_n$ has the form

(A.22)
$$W_n = \begin{pmatrix} 0 & 0 & \cdot & \cdot & \cdot & 1 \\ 0 & 0 & \cdot & \cdot & 1 & 0 \\ 0 & 0 & \cdot & \cdot & 0 & \cdot \\ 0 & 0 & \cdot & \cdot & 0 & \cdot \\ 0 & 1 & \cdot & \cdot & 0 & 0 \\ 1 & 0 & \cdot & \cdot & \cdot & 0 \end{pmatrix}.$$

Note that the matrix $C_n(f) - T_n(f)$, the difference between a Toeplitz matrix and its associated circulant, is a $90°$ clockwise rotation of the matrix $H_n(f) + W_n H_n(f) W_n$, where $H_n(f)$ is the Hankel matrix with symbol $f$ and $W_n$ is the flip matrix as given in (A.22). Part (b) follows from the fact that all the elements of $H_n(f) + W_n H_n(f) W_n$ are zero except for the first and the last principal submatrix of size $N \times N$.

Now the singular values of $C_n(f) - T_n(f)$ and $H_n(f) + W_n H_n(f) W_n$ are the same. Since $W_n$ is an orthogonal matrix, the singular values of $H_n(f)$ and $W_n H_n(f) W_n$ are the same. Now, part (c) follows from an application of part (b) of Lemma 1. $\square$

PROOF OF LEMMA 4. First note that element $(j, k)$ of the matrix $T_n(s^d)$ is given by $(-1)^{j-k} \binom{2d}{d-|j-k|}$. It is then enough to show that, for any $d+1 \leq j, k \leq n - d$, element $(j, k)$ of the matrix $\overline{S}'_{-d} \overline{S}_{-d}$ is $(-1)^{j-k} \binom{2d}{d-|j-k|}$.

It is not difficult to see that the following identity is valid:

(A.23)
$$\sum_{0 \leq t \leq l} \binom{d}{t} \binom{d}{l-t} = \binom{2d}{l}.$$

This identity follows from expanding $(1-z)^{2d}$ as $\sum_{0 \leq t \leq 2d} (-1)^t \binom{d}{t} z^t$. On the other hand, we can expand $(1-z)^{2d}$ as

$$(1-z)^d (1-z)^d = \left\{ \sum_{0 \leq s \leq d} (-1)^s \binom{d}{s} z^s \right\} \left\{ \sum_{0 \leq t \leq d} (-1)^t \binom{d}{t} z^t \right\}.$$

Note that element $(j, k)$ of the matrix $\overline{S}'_{-d} \overline{S}_{-d}$ is given by

$$\sum_{1 \leq t \leq n-d} (-1)^{t-j} \binom{d}{t-j} (-1)^{t-k} \binom{d}{t-k}$$

$$= (-1)^{j-k} \sum_{1 \leq t \leq n-d} \binom{d}{t-j} \binom{d}{t-k}.$$

Now, use of the identity (A.23) on the right-hand side of the last expression yields the desired result. $\square$

**Acknowledgments.** The authors are indebted to the Associate Editor and three referees who made valuable comments and suggestions. These lead us to investigate some issues not addressed in the first draft of the paper, and the revision contains answers to these issues.

Department of Statistics
University of California
Davis, California 95616
USA
E-mail: burman@wald.ucdavis.edu
shumway@wald.ucdavis.edu